\documentclass[12pt,psamsfonts]{article}
\input epsf
\usepackage{psfig}
\usepackage{latexsym}
\usepackage{amsfonts}
\usepackage[psamsfonts]{amssymb}
\def\dim{{\mathrm{dim}\,}}
\def\codim{{\mathrm{codim}\,}}
\def\rank{{\mathrm{rank}\,}}
\def\card{{\mathrm{card}\,}}

\def\s{{\mathbf{s}}}
\def\x{{\mathbf{x}}}

\def\T{{\mathbf{T}}}

\def\a{\mbox{\boldmath$a$}}
\def\scripta{\mbox{\boldmath$\scriptstyle a$}}

\def\CP{{\mathbf{P}}}
\def\RP{{\mathbf{RP}}}
\newcommand{\binom}[2]{\left(\begin{array}{c}#1\\#2\end{array}\right)}

\title{Rational functions and real Schubert calculus}
\author{
A.\ Eremenko\thanks{
supported by NSF grants DMS-0100512 and DMS-0244421.},
A.\ Gabrielov\thanks{Supported by NSF grants DMS-0200861 and
DMS-0245628.},
M.\ Shapiro\thanks{Supported by NSF grant 0401178, BSF grant 2002375 and by
the Institute of Quantum Science, MSU.} and
A.\ Vainshtein\thanks{Supported by BSF grant 2002375.}}
\date{\today}
\begin{document}
\maketitle

\noindent
{\bf 1. Introduction}
\vspace{.1in}

We single out some problems of Schubert calculus of subspaces
of codimension $2$
that have the property that all their solutions are real whenever the
data are real. 
Our arguments explore the connection between subspaces of codimension $2$
and rational functions of one variable.

Let $\{ W_j\}_{j=1}^q$ be a finite collection of projective subspaces
in general position in the complex
projective space $\CP^d,\; a_j=\dim W_j$,
and
\begin{equation}
\label{1}
1\leq a_j\leq d-1,\quad \sum_{j=1}^qa_j=2d-2.
\end{equation}
We consider the following
\vspace{.1in}

\noindent
{\bf Problem 1}.  {\em Enumerate subspaces $X\subset\CP^d$ of codimension $2$ which
intersect each $W_j$ non-transversally, that is}
\begin{equation}
\label{2}
\dim X\cap W_j\geq a_j-1\quad\mbox{\em for every}\quad j\in[1,q].
\end{equation}
If the given subspaces $W_j$ are in general position then the number of such
complex subspaces $X$ is finite and can be obtained by using Pieri's formula
from Schubert calculus \cite[Ch. I, sect. 5]{gh}.

This number turns out to be the {\em Kostka number} corresponding
to the shape $2\times(d-1)$
whose definition we recall. Let $\a=(a_1,\ldots,a_q)$.
Consider the Young diagrams of the shape $2\times(d-1)$.
They consist of two rows of length $d-1$.
A {\em semi-standard Young tableau} SSYT of shape $2\times(d-1)$
is a filling of such a
diagram by positive integers, such that an integer $k$ appears
$a_k$ times, the entries are strictly increasing in the columns and
non-decreasing in the rows.
The corresponding Kostka number
$K_{\scripta}$
is the number of such SSYT.
Schubert calculus interpretation of these numbers shows that
$K_{\scripta}$ does not change if we permute the coordinates of $\a$.
For a purely combinatorial proof of this see \cite[Cor. 1.2.9]{ma}
or \cite[Theorem 7.10.2]{stan}.

In this paper we treat Problem~1 over the field of
real numbers. A subspace $X\subset\CP^d$ is called
{\em real} if it can be defined by equations with real
coefficients, or, equivalently, if $X$ is generated by a set
of real points. It follows from  a general result of Sottile \cite{sot}
that {\em there exist} configurations of real subspaces $W_j$
such that all solutions $X$ of Problem~1 are real.
The question is to find an explicit condition on $W_j$ which would
imply that all solutions are real.

The following condition was proposed by B. and M. Shapiro
(see \cite{shap1,shap2}).
Let $E:\CP^1\to\CP^d$ denote the rational normal curve; in homogeneous
coordinates
$E(z)=(1:z:\ldots:z^d)$. The B. and M. Shapiro conjecture for the case
of subspaces of codimension $2$
says that
if all $W_j$ are osculating $E$ at distinct real points then all
solutions $X$ of Problem~1 are real. This conjecture
was proved in \cite{eg}. B. and M. Shapiro made a similar conjecture
for enumerative problems involving subspaces of arbitrary codimension
but in this paper we only consider the case of codimension $2$.

In this paper we consider an extension
of the result in \cite{eg}
to the case that $W_j$ are {\em spanned} by some finite
sets of real points on $E$.

Let $\{ A_j\}_{j=1}^q$ be a collection of finite sets on the circle
$\RP^1\subset\CP^1$.
We say that this collection is {\em separated} if there
exist pairwise disjoint closed arcs $I_j\subset\RP^1$
such that $A_j\subset I_j$
for $1\leq j\leq q$. Our main result is

\vspace{.1in}
\noindent
{\bf Theorem 1} {\em  Let $\{ A_j\}_{j=1}^q$ be a separated collection
of finite sets in $\RP^1$, each set $A_j$ containing $a_j+1$
points, such that $(\ref{1})$ holds,
and let $W_j$ be subspaces in $\CP^d$ spanned by the sets $E(A_j)$.
Then all subspaces $X$ of codimension $2$ which satisfy $(\ref{2})$ are real.
For a generic configuration of $A_j$ there are exactly $K_{\scripta}$ such
subspaces $X$.}
\vspace{.1in}

In general, the condition that $A_j$ are separated cannot
be removed. Let us consider several special cases.
\vspace{.1in}

\noindent
1. Let $q=2$, $a_1=a_2=d-1$. Then the problem always has one solution,
this solution is real, and the condition that $A_j$ are separated is
redundant.
\vspace{.1in}

\noindent
2. Let us consider the limiting
situation when all points in each $A_j$ collide.
Here $a_j=\card A_j-1$
are arbitrary integers satisfying (\ref{1}).
Then $W_j$ are subspaces osculating $E$
at some real points;
this case was considered in \cite{eg,eg2}.
We will show that this situation is generic enough and
the number of solutions is exactly $K_{\scripta}$, for {\em any}
choice of the real points (Theorem 3 below).
If $\a=(1,\ldots,1)$ then $K_{\scripta}$ is the Catalan number
\begin{equation}
\label{cat}
K_{(1,\ldots,1)}=\frac{1}{d}\binom{2d-2}{d-1}.
\end{equation}
\vspace{.1in}

\noindent
3. Now we consider the case that $a_1\in [1,d-1]$ is arbitrary
and all the rest of $a_j$ are equal to $1$.
The Kostka number in this case is
\begin{equation}
\label{num}
K_{(a_1,1,\ldots,1)}=\frac{a_1+1}{d}\binom{2d-2-a_1}{d-1}.
\end{equation}
In the case that the two points of each $A_j,\; j\geq 2$ collide,
Problem~1 is equivalent to the following
problem of enumeration of flags.

Consider flags $F=(F_1,F_2)$ in $\CP^d$, where
$$F_2\subset F_1\quad\mbox{and}\quad\codim F_i=i,\quad i=1,2.$$
Suppose that $A_1=\{ x_0,\ldots,x_{a_1}\}\subset\RP^1$ and $A_j=\{ y_j^{(2)}\}$,
where the superscript $2$ indicates that each
point $y_j$ is taken with multiplicity $2$, $j=2,\ldots, 2d-2-a_1$.
Our condition on the flags $F$ is that $F_1$ contains the points
$E(x_k),\; 0\leq k\leq a_1$, and $F_2$ intersects the tangent lines $W_j$
to $E$ at $y_j$, $j=2,\ldots, 2d-2-a_1$.
We want to enumerate such flags $F$.

The condition on $F_1$ implies that $F_1$ contains $W_1$,
the subspace spanned by $E(x_k)$, $0\leq k\leq a_1$.
If $x_k$ are in general position
then $\dim W_1=a_1$. Now $W_1\subset F_1$ and $F_2\subset F_1$ imply
that
\begin{equation}
\label{4}
\dim W_1\cap F_2\geq a_1-1,
\end{equation}
and the remaining conditions on $F_2$ are
\begin{equation}
\label{5}
\dim F_2\cap W_j\geq a_j-1=0,\quad j=2,\ldots,2d-2-a_1.
\end{equation}
In the opposite direction, suppose that $F_2$ satisfies (\ref{4}) and (\ref{5}),
and define $F_1$ as the span of $F_2$ and $W_1$.
Then $\codim F_1=1$ and $F_2\subset F_1$.

Thus our problem of flag enumeration is equivalent to the
problem of enumeration of subspaces $F_2$ of codimension $2$ satisfying
(\ref{4}) and (\ref{5}), which is a special limiting case of Problem 1.
Assuming that the sets $A_j$
are separated, Theorem~1 implies that all solutions are real,
and for generic choice of points $x_i$ and $y_j$ the number of these
solutions is given by (\ref{num}).
This problem of enumeration of flags was subject to extensive
computer experiments \cite{sot7}.

In these experiments the role of the separation condition was discovered:
when it holds all solutions are real, while when it does not hold,
many configurations of points $A_j$ give fewer real solutions than
the upper estimate from Schubert calculus.

All authors thank Frank Sottile for stimulating discussions
and MSRI for the opportunity to work together during the
semester ``Topological methods in real algebraic geometry'' in
spring 2004.
\vspace{.2in}

\noindent
{\bf 2. Rational functions}
\vspace{.1in}

The method of this paper is based on the relation
between subspaces of codimension $2$ and rational functions
of one variable.

Let $G=G(d-1,d+1)$ be the Grassmannian of projective
subspaces of dimension $d-2$ in $\CP^d$. Every such subspace $X$
can be defined by two equations in homogeneous coordinates
\begin{equation}
\label{X}\begin{array}{rcl}b_{0,0}z_0+\ldots+b_{0,d}z_d&=&0,\\
                  b_{1,0}z_0+\ldots+b_{1,d}z_d&=&0.
\end{array}
\end{equation}
We put into correspondence to this system a rational function
\begin{equation}
\label{xx}
f(z)= \frac{b_{0,0}+b_{0,1}z+\ldots+b_{0,d}z^d}{b_{1,0}+b_{1,1}z+\ldots+
b_{1,d}z^d}.
\end{equation}
Choosing another system of equations of the form
(\ref{X}) that defines the same $X$
results in
replacing
$f$ by
$\phi\circ f$, where $\phi$ is a rational function of degree $1$.
We will call rational functions $f_1$ and $f_2$ {\em equivalent}
if $f_1=\phi\circ f_2,\;\deg\phi=1$. Let $D\subset G$ be the
subset that corresponds to reducible rational functions
(those whose numerator and denominator have a non-constant
common factor). Then $G\backslash D$ is in bijective correspondence
with equivalence classes of non-constant rational functions of degree at
most $d$.

A rational function is called {\em real} if it maps $\RP^1$ into
itself. An equivalence class is called real if it contains a real
rational function.
Real classes correspond to the real elements of $G$.

Now we translate conditions (\ref{2}) to the language of rational
functions. Let $W$ be the subspace of dimension $a$ spanned
by the points $E(z_k), \; 0\leq k\leq a$ on the rational normal curve.
The condition that
$X$ given by (\ref{X}) intersects $W$ non-transversally means
that
\begin{equation}
\label{eq}
\rank BC=1,
\end{equation}
where
\begin{equation}
\label{mat}
B=\left(\begin{array}{lll}b_{0,0}&\ldots&b_{0,d}\\
                          b_{1,0}&\ldots&b_{1,d}
\end{array}\right),
\end{equation}
and $C$ is the matrix whose columns are $(1,z_k,\ldots,z_k^d),$
$0\leq k\leq a$. Equation (\ref{eq}) says that the two
rows of $BC$ are proportional, which is equivalent to
\begin{equation}
f(z_0)=f(z_1)=\ldots=f(z_{a}).
\end{equation}
Thus in the case that each $W_j$ is the span of $a_j+1$ points on
the rational normal curve, Problem~1
is equivalent to the following interpolation-type
problem for rational functions:
\vspace{.1in}

\noindent
{\bf Problem} $\mathbf{1^\prime}$. {\em For given subsets $A_j\subset\CP^1,$ of
cardinality $a_j+1$, where $a_j$ satisfy $(\ref{1})$,
enumerate classes of rational functions $f$ of degree $d$
with the property that for every $j$, $f$ is constant on $A_j$.}
\vspace{.1in}

Thus our Theorem~1 is equivalent to
\vspace{.1in}

\noindent
{\bf Theorem $\mathbf{1^\prime}$} {\em If the sets $A_j\subset\RP^1$
are separated then
all solutions of Problem~$1'$ are real.
For generic subsets $A_j\subset\RP^1$ the number of classes of real
rational functions
solving Problem~$1'$ is $K_{\scripta}$.}
\vspace{.1in}

It is not known whether the genericity condition
in the second statement of Theorem~${1^\prime}$ can be removed.
It can be removed in the following special case.
\vspace{.1in}

\noindent
{\bf Theorem 2} {\em Let $A_1,\ldots,A_q$ be separated sets on the real line,
$\card A_j=a_j+1\in [2,d]$, and
$$\sum_{j=1}^qa_j=d-1.$$
Then there exists a unique class of polynomials $p$ of degree $d$,
satisfying the conditions that $p$ is constant on $A_j$
for $1\leq j\leq q$.}
\vspace{.1in}

{\em Proof.} We restrict ourselves
 for simplicity to the case that $q=d-1$ and
$a_j=1$ for
$1\leq j\leq q$. The general case is similar. Let $A_j=\{ x_j,y_j\}$,
where
\begin{equation}
\label{poly2}
x_1<y_1<x_2<y_2<\ldots<x_{d-1}<y_{d-1}.
\end{equation}
Every polynomial of degree at most $d$ is equivalent
to a unique polynomial of the form
\begin{equation}\label{pol1}
p(z)=z^n+\ldots+b_1z,\quad\mbox{where}\quad n\leq d.
\end{equation}
Conditions of the theorem mean that
\begin{equation}
\label{poly1}
p(x_j)=p(y_j),
\end{equation}
for $1\leq j\leq d-1$. This is a system of linear equations
with respect to coefficients of $p$, and it is easy to see that
the coefficients of this system are real.
Thus the system (\ref{poly1}) has unique solution of the form
(\ref{pol1}) if and only if it has a unique real solution of this form.

Now suppose that our system (\ref{poly1}) has two solutions of the
form (\ref{pol1}). Subtracting these solutions we obtain another solution
which is a real polynomial $p_0$ of degree strictly less than $d$,
and
\begin{equation}
\label{noll}
p_0(0)=0.
\end{equation}
On the other hand, Rolle's theorem and (\ref{poly1}) implies that
the derivative $p^\prime_0$ has at least $d-1$ zeros.
So $p^\prime_0=0$ and we conclude from (\ref{noll}) that $p_0=0$.
This proves the theorem.
\hfill $\Box$
\vspace{.1in}

In the limit, when in each group $A_j$ all points collide to
one point $x_j$ we recover from Theorem $1'$ the main result of \cite{eg}:
{\em all rational functions whose critical points are real are equivalent
to real rational functions.}

In this degenerate case we have an important additional information:
\vspace{.1in}

\noindent
{\bf Theorem 3} {\em
The number of classes of rational functions of degree $d$ with
arbitrarily prescribed real critical points of
given multiplicity is exactly $K_{\scripta}$, for any choice of the real
critical points,
where $\a=(a_1,\ldots,a_q)$ is the vector of multiplicities.}
\vspace{.1in}

Our proofs of Theorems $1'$ and 3 are
based on the results in \cite{eg,eg2} which we recall in the next section.
\vspace{.2in }
\newpage

\noindent
{\bf 3. Proof of Theorems $\mathbf{1^\prime}$ and 3}
\vspace{.1in}

It is sometimes convenient to replace
$\RP^1$ by the unit circle $\T$ as in \cite{eg}.
We always assume that $\T$ is equipped with the counter-clockwise orientation.
We denote by $\s:\CP^1\to\CP^1$ the reflection with respect to $\T$,
that is $\s(z)=1/\overline{z}$.

We fix $d\geq 2$.
Let $R$ be the class of all rational functions $f$ such that $f(\T)\subset\T$,
all critical points of $f$ are simple and belong to $\T$. We equip $R$
with the following topology (see \cite{eg2}):
a sequence $(f_k)$ converges to $f$ if there exists a finite
set $S\subset\CP^1$ such that $f_k\to f$ uniformly on compact subsets
of $\CP^1\backslash S$.

First we recall a parametrization of the class $R$. Let $f$ be a function
of the class $R$, and $v_0$ a critical point of $f$.
The full preimage $\Gamma=f^{-1}(\T)$ defines
a cell decomposition of $\CP^1$. We
call the pair $\gamma=(\Gamma,v_0)$ the {\em net of $f$ with respect
to $v_0$.}

One can describe all objects which can occur in this way.
Let $\Gamma$ be a cell decomposition of $\CP^1$
with the following properties:
\vspace{.1in}

\noindent
(i) closure of each cell is homeomorphic to a closed ball
of corresponding dimension,
\vspace{.1in}

\noindent
(ii) $\Gamma$ is symmetric with respect to $\T$, that is the 
reflection $\s$ maps each cell of $\Gamma$ onto a cell of $\Gamma$,
\vspace{.1in}

\noindent
(iii) all $0$-cells belong to the unit circle $\T$ and there
are $2d-2$ of them,
\vspace{.1in}

\noindent
(iv) $1$-skeleton of $\Gamma$ contains $T$.
\vspace{.1in}

As usual we call $0$-, $1$- and $2$ cells vertices, edges and
faces. The edges that do not belong to $\T$ will be called
{\em interior edges}.
A {\em net} is a pair $\gamma=(\Gamma,v_0)$ where
$\Gamma$ is a cell decomposition of $\CP^1$ with the properties
(i)-(iv) and $v_0$ is a vertex of $\Gamma$.

Two nets $(\Gamma,v_0)$ and $(\gamma',v_0^\prime)$
are called {\em equivalent} if there exists a homeomorphism
$\phi$ of $\CP^1$ preserving the orientations of both
$\CP^1$ and $\T$, commuting with the reflection $\s$
and having the properties $\phi(\Gamma)=\Gamma'$
and $\phi(v_0)=v_0^\prime$.

One of the results of \cite{eg} is the following: {\em Let
$f_1$ and $f_2$ be  two rational
functions in $R$ with the same critical points, and $v_0$
one of these common critical points. If the nets
of $f_1$ and $f_2$ with respect to $v_0$ are equivalent then
$f_1$ and $f_2$
are
equivalent.}

The main technical result of \cite{eg} can be stated as follows:
{\em For every net $\gamma$ whose vertex set is $V$,
and for every injective map $c:V\to\T$
preserving the cyclic order induced by the orientation, 
there exists a unique class $C$ of rational functions in $R$
whose critical points are $c(V)$
and whose net with respect to $c(v_0)$ is equivalent to $\gamma$.}

We will say somewhat informally that functions of this class
$C$ have prescribed net and prescribed critical points.

Functions $f\in R$ with a prescribed net
and prescribed critical points
can be normalized in the following way. Choose two additional
vertices $v_1$ and $v_{-1}$ and three distinct points
$w_{-1},w_0,w_1$ on $\T$. Then the normalization condition is
\begin{equation}
\label{norm}
f(c(v_i))=w_i,\quad i=-1,0,1.
\end{equation}
There exists a unique normalized function $f\in R$ with a
prescribed net and prescribed critical points. Another result from
\cite{eg} says that for each fixed net, this normalized function
depends continuously on the critical points, or more precisely,
on the map $c$ above.

Now we can prove Theorem~$1'$.
Suppose that the separated sets $A_j\subset\T$
satisfying (\ref{1}) are given.
Let $I_j\subset\T$ be disjoint arcs containing $A_j$. We list each $A_j$ as
$(z_{j,0},z_{j,1},\ldots,z_{j,a_j})$ where the order
within $A_j$ is consistent
with the counter-clockwise orientation of $\T$. For any two distinct
points $r$ and $s$ in $\T$,
we denote by $[r,s]$ the
closed arc of $\T$ which is traced from $r$ to $s$
counter-clockwise.

Choose additional points
\begin{equation}
\label{points}
x_{j,k}\in[z_{j,k},z_{j,k+1}],\quad 0\leq k\leq a_j-1,\quad 1\leq j\leq q.
\end{equation}
According to (\ref{1}), the total number of these points $x_{j,k}$
is $2d-2$. So we can use these points as vertices of a net.
We consider the nets satisfying the following additional condition:
\begin{equation}
\label{cond}
\mbox{\it there are no interior edges connecting two vertices on
the same $I_j$.}
\end{equation}

{\bf Lemma 1} {\em Suppose that a net $\gamma$ satisfies $(\ref{cond})$.
Let all $x_{j,k}$, except one, $x=x_{j_0,k_0}$
be fixed arbitrarily, so that
conditions $(\ref{points})$ are satisfied.
For fixed $w_i\in\T$, $j=0,\pm1$,
let $f_x$ be the function of the class
$R$ with the net $\gamma$ and critical points $x_{j,k}$,
normalized as in $(\ref{norm})$. Denote $I=[z_{j_0,k_0},z_{j_0,k_0+1}].$
Then there exist a continuous branch of $\arg f_x(z)$,
for $(x,z)\in I\times I$, such that the function
$$\psi(x)=\arg f_x(z_{j_0,k_0})-\arg f_x(z_{j_0,k_0+1})$$
changes sign as $x$ runs over $I$.}

{\em Proof.} To simplify our notation we put $r=z_{j_0,k_0},
s=z_{j_0,k_0+1}$, so that $I=[r,s]$. Let $\ell$ be the interior edge
of $\gamma$ that lies in the unit disc and has $x$ as an endpoint.
Let $v$ be another endpoint of this edge. We note that
neither the assumptions nor the conclusions of the lemma depend on
the choice of the normalization. So we may assume that $v$ is one
of the distinguished vertices, say $v=v_1$ in (\ref{norm}) and
choose $w_1=f_x(v)=-1.$
Let $D_1$ and $D_2$ be two faces of the net $\gamma$ in the unit disc which have
$\ell$ as the common part of their boundaries. As the boundary of every face
is mapped by $f_x$ bijectively onto $\T$,
and the vertex $v$ is disjoint from $I$ because of (\ref{cond})
we conclude that $f_x(z)\neq -1$
for $(x,z)\in I\times I.$
Thus
$\arg f_x(z)$ can be defined as a continuous function for
$(x,z)\in I\times I$,
with values in $(-\pi,\pi)$.
Now, $z\mapsto \arg f_x(z)$ is strictly monotone in the opposite directions
on the arcs $[r,x]$ and $[x,s]$: it is strictly increasing
on one of these arcs and strictly decreasing on another one.
This follows from the fact that one of the closed faces $D_1$ or $D_2$
is mapped homeomorphically onto the closed unit disc and another
onto the closed exterior of the unit disc, while the arcs $[r,x]$ and $[x,s]$
are parts of the boundaries of $D_1$ and $D_2$.

We conclude that $\arg f_x(r)-\arg f_x(s)$ changes sign as $x$ runs from
$r$ to $s$, which proves the lemma.\hfill$\Box$
\vspace{.1in}

\noindent
{\bf Lemma 2} {\em If $\gamma$ is a net satisfying
$(\ref{cond})$ then
there exists a choice of critical points $x_{j,k}$ such that every
function $f\in R$ with this net and with these critical points
satisfies}
\begin{equation}
\label{ne}
f(z_{j,0})=\ldots=f(z_{j,a_j})\quad\mbox{for every}\quad j\in [1,q].
\end{equation}

{\em Proof}. Let $\x$ be the vector with coordinates $(x_{j,k})$.
The set of $\x$ satisfying (\ref{points}) is a closed cube $Q$.
According to Lemma 1,
the continuous functions
$\psi_{j,k}(\x)=\arg f_{\x}(z_{j,k})-\arg f_{\x}(z_{j,k+1})$
with appropriate choices of branches of arguments take values
of opposite signs on the opposite facets of $Q$.
It follows from a Corollary of Brouwer's Theorem \cite[Ch. IV, 1D]{gw}
that all these functions $\psi_{j,k}$ have
a common zero in the interior of $Q$.
\hfill$\Box$
\vspace{.1in}

So, for each net $\gamma$ satisfying (\ref{cond}),
Lemma~2 gives at least one function $f\in R$ that satisfies
(\ref{ne}). The functions with different nets cannot be equivalent.
The crucial fact is
\vspace{.1in}

\noindent
{\bf Lemma 3} {\em For each choice of the separated sets $A_j$ with
$\card A_j=a_j+1$
there are exactly $K_{\scripta}$ nets satisfying $(\ref{cond})$.}
\vspace{.1in}

{\em Proof}. We describe two injective maps:
one from the set of nets to the set of SSYT (part A below) and
another one in the opposite direction (part B). 

First of all we observe that
if two nets $(\Gamma,v_0)$ and $(\gamma',v_0^\prime)$ satisfying
(\ref{cond}) are equivalent then $v_0$ and $v_0^\prime$
cannot belong to the same interval $I_j$.
Therefore, instead of distinguishing a vertex $v_0$, it is sufficient
to distinguish one of the intervals $I_j$. This is achieved
by considering the real line instead of the unit circle $\T$,
which distinguishes the leftmost interval.

We assume for convenience that
the intervals are enumerated in such a way
that $i<j$ implies $x<y$ for all $x\in I_i$ and $y\in I_j$,
besides the vertices in each interval are enumerated from left to
right.
\vspace{.1in}

\noindent
{\em Part A.}
We construct an SSYT corresponding to a given net. 
To each vertex of the net we assign one entry in SSYT.
These entries are defined inductively: for $n\geq 1$, suppose that
the entries corresponding to $x_k$ for $k\leq n-1$ are already
placed in the tableau. The entry corresponding to $x_n$
is the number $j$ such that $x_n\in I_j$. This number $j$
is placed in the first row if the edge from $x_n$ goes forward
(to some $x_m>x_n$) and in the second row if the edge from $x_n$
goes backward. The new entry $j$ is placed to the leftmost place
in its row. It is easy to see that in this way we obtain an SSYT
for each net satisfying (\ref{cond}).
\vspace{.1in}

\noindent
{\em Part B.}
In the opposite direction, we construct inductively the net 
corresponding to a given SSYT of shape $2\times(d-1)$.
On the $n$-th step of the construction we deal with a tableau of
shape $2\times(d-n)$ filled by positive integers such that the rows
are non-decreasing and columns are strictly increasing.
The only difference from SSYT is that the integers in the tableau
are not necessarily consequtive; slightly abusing terminology
we will call these tableaux SSYT as well.
Assuming that such an SSYT is given, let $k$ be the first
(leftmost) entry in its second row.
We find the rightmost entry $m$ in the first row that
is strictly less than $k$ and draw an edge between
the leftmost free vertex
in $I_k$ and the rightmost free vertex of $I_m$.
Then we delete the first entry from the second row and the
above defined entry from the first row, and shift the left part
of the first row to the right to form a tableau of
shape $2\times(d-n-1)$.
Evidently, the obtained tableau is an SSYT, so we can proceed until
it becomes empty.

To complete Part B, it remains to prove that
the edge added on the $n$-th step does not
intersect the edges added on the previous steps.

Let us show first that on each step all the vertices lying
between the endpoints  
of the added edge are not free. Indeed, suppose that there is a free
vertex $v$ between the endpoints. We use the notation
from the description of the $n$-th step above.
By assumption, the vertex $v$ belongs to some $I_p$
with $m<p<k$, since the new edge joins the leftmost free vertex
in $I_k$ to the rightmost free vertex in $I_m$.
Since the vertex $v$ in $I_p$ is free there should be an entry $p$
somewhere in the current SSYT. It cannot be in the second row since 
$k>p$ is the leftmost entry of the second row.
It cannot be in the first row since $m<p$ is the rightmost entry
of those which are strictly less than $k$. 

Now we can complete the proof. If a new edge intersects am old one,
then one of the vertices of the new edge was a free vertex
between the endpoints of the old edge on the step when this old
edge was added.
This completes the proof of the part B. 
\hfill $\Box$
\vspace{.1in}

Now we can complete the proof of Theorem $1'$.
For given sets $A_j$ in general position, Schubert calculus gives
at most $K_{\scripta}$ classes of rational functions solving Problem $1'$.
On the other hand, Lemmas 2 and 3 give at least $K_{\scripta}$
real classes of such functions. Thus all classes are real.
Passing to the limit, we obtain that for all separated sets $A_j$,
all solutions of Problem $1'$ are real.

\hfill$\Box$
\vspace{.1in}

{\em Proof of Theorem 3.} Let $R'$ be the subset of the closure of $R$
which consists of functions of degree exactly $d$.
The nets of functions $f\in R'$ have similar properties to those of
$f\in R$, except that the vertex degrees can be arbitrary numbers
between $4$ and $2d$. Such nets were called degenerate in \cite{eg2}.
In that paper, the following result was proved:
{\em Let $x_1^{a_1},\ldots,x_q^{a_q}$ be a divisor on the real
line, where $\a=(a_1,\ldots,a_q)$ satisfies $(\ref{1})$
and $\gamma$ be a net whose vertices enumerated in the increasing order
have degrees $2a_j+2$. Then there exists a real rational function of degree
$d$ with critical points of multiplicity $a_k$ at the points $x_k$.}

It is easy to see that the number of these degenerate nets $\gamma$ is the
same as the number of non-degenerate nets satisfying (\ref{cond}),
so according to Lemma 3, this number is $K_{\scripta}$. Hence
the number of classes
of real functions satisfying the conditions of Theorem 3 is at least
$K_{\scripta}$ and this proves the theorem. \hfill$\Box$

{\em A. E.: Purdue University, West Lafayette, IN 47907-2067.

eremenko@math.purdue.edu

A. G.: Purdue University, West Lafayette, IN 47907-2067.

agabriel@math.purdue.edu

M. S.: Michigan State University, East Lansing, MI 48824.

mshapiro@math.msu.edu

A. V.:University of Haifa, Mount Carmel, 31905 Haifa, Israel.

alek@cs.haifa.ac.il}

\vspace{.2in}

\begin{thebibliography}{1}

\bibitem{eg} A. Eremenko and A. Gabrielov,
Rational functions with real critical points
and the B. and M. Shapiro conjecture in real enumerative geometry,
Ann. Math., 155 (2002), 105-129.
\bibitem{eg2} A. Eremenko and A. Gabrielov,
Wronski map and Grassmannians of real codimension 2 subspaces,
Computational Methods and Function Theory, 1 (2001) 1-25.
\bibitem{gh} Ph. Griffiths and J. Harris, Principles of algebraic
geometry, Willey-Interscience, NY, 1978.
\bibitem{gw} W. Hurewicz and H. Wallman, Dimension theory, Princeton
Univ. Press, 1948.
\bibitem{ma} L. Manivel, Symmetric functions, Schubert Polynomials and
degeneracy loci, AMS, Providence, RI, 2001.
\bibitem{sot7} J. Ruffo, Y. Sivan, E. Soprounova and F. Sottile,
Experimentation and conjectures in the real Schubert calculus,
work in progress,
www.math.umass.edu/{\~{ }}sottile/pages/Flags
\bibitem{sot} F. Sottile, The special Schubert calculus is real,
Electron. Res. Announc. AMS, 5 (1999) 35--39.
\bibitem{shap1} F. Sottile, Real Schubert calculus: polynomial systems and a
conjecture of Shapiro and Shapiro, Experimental math., 9 (2000) 161--182.
\bibitem{shap2} F. Sottile, Enumerative real algebraic geometry,
DIMACS Series in Discrete math. and Computer sci., vol. 60,
AMS, Providence, RI, (2003) 139--179.
\bibitem{stan} R. Stanley, Enumerative Combinatorics, vol. 2,
Cambridge UP, 1999.
\end{thebibliography}
\end{document}